\documentclass{amsart}
\usepackage{ryanmacros}
\input{ryanquivers.sty}
\usepackage{fullpage}

\usepackage[colorlinks=true, pdfstartview=FitV, linkcolor=blue, citecolor=blue, urlcolor=blue]{hyperref}

\CompileMatrices

\begin{document}

\title{The Rank of a Quiver Representation}
\author[Kinser]{Ryan Kinser}
\address{Department of Mathematics, University of Michigan, Ann Arbor, Michigan 48109}
\email{kinser@umich.edu}
\thanks{This material was based upon work supported under a National Science Foundation Graduate Research Fellowship, and NSF Grant DMS 0349019}

\begin{abstract}
We define a functor which gives the ``global rank of a quiver representation'' and prove that it has nice properties which make it a generalization of the rank of a linear map.  We demonstrate how to construct other ``rank functors'' for a quiver $Q$, which induce ring homomorphisms (called ``rank functions'') from the representation ring of $Q$ to $\Z$.  These rank functions give discrete numerical invariants of quiver representations, useful for computing tensor product multiplicities of representations and determining some structure of the representation ring.  We also show that in characteristic 0, rank functors commute with the Schur operations on quiver representations, and the homomorphisms induced by rank functors are $\lambda$-ring homomorphisms.  
\end{abstract}
\maketitle

\comment{
introsect
defssect
repringsect
constructionsect
ranktensorsect
}



\section{Introduction}\label{introsect}

Fix a field $K$.  For a quiver $Q$, let $\repq = \rep_K(Q)$ be the category of finite dimensional representations of $Q$ over $K$.  We will take all quivers to be connected throughout.  There is a natural tensor product of quiver representations (Section \ref{tensorsect}), giving $\repq$ the structure of a tensor category.  If $Q$ and $Q'$ are quivers, we will call a functor $F: \repq \to \rep (Q')$ a \textbf{tensor functor} (Section \ref{categorysect}) if it satisfies:
\begin{enumerate}[(a)]
\item For representations $V$ and $W$, there are isomorphisms
\[
F(V) \oplus F(W) \cong F(V \oplus W) \qquad \text{and} \qquad F(V) \otimes F(W) \cong F (V \otimes W)
\]
which are functorial in $V$ and $W$, and symmetric.
\item $F(\id_Q) \cong \id_{Q'}$, where $\id_Q$ is the identity representation of $Q$ (Section \ref{tensorsect}).
\end{enumerate}

For example, let $p$ be a path in some quiver $Q$.  The quiver $A_1$ has a single vertex and no arrows, so  $\rep(A_1) = \vscat$ is the category of finite dimensional $K$-vector spaces.  Then the functor
\[
\imfunc_p : \repq \to \rep(A_1)		\qquad V \mapsto \im V_p 
\]
is a tensor functor. In particular, the trivial path $\varepsilon_x$ at a vertex $x$ corresponds to the functor $\imfunc_{\varepsilon_x} (V) = V_x$.  But in general, these are not the only tensor functors on $Q$.  Motivated by this example, we will call a tensor functor from $\repq$ to $\rep(A_1)$ a \textbf{rank functor} on $Q$.

For a fixed quiver $Q$, the operations $\oplus$ and $\otimes$ give a semiring structure to the set of isomorphism classes of representations of $Q$.  By introducing virtual isomorphism classes of representations, we can form a ring $R(Q)$ from this semiring, called the \textbf{representation ring} of $Q$ (Section \ref{repringsect}).  Such a technique was introduced by Grothendieck (communicated in \cite{BSriemannroch}) for coherent sheaves on an algebraic variety, and has since appeared, for example, in equivariant K-theory \cite{segalequivktheory} and to study representations of Lie groups \cite{segallierepring}.
When $\charac K =0$, we can define Schur functors on $\repq$, giving $R(Q)$ the structure of a $\lambda$-ring \cite{knutsonlambda}.  In the quiver setting, the structure of $R(Q)$ has been determined by Herschend for $Q$ of type $A_n$ and $D_n$ in \cite{herschend07b}, and for type $\widetilde{A_n}$ \cite{herschendaffinean}.

A tensor functor $F: \repq \to \rep(Q')$ induces a ring homomorphism $f: R(Q) \to R(Q')$, which is a homomorphism of $\lambda$-rings in characteristic 0 (Theorem \ref{rankschurthm}).  In particular, rank functors on $Q$ induce ring homomorphisms $R(Q) \to \Z$ via the isomorphism $R(A_1) \cong \Z$ which identifies a vector space with its dimension.  These homomorphisms give numerical invariants of quiver representations, and can be used to study the structure of $R(Q)$.

Given any quiver $Q$, we will construct in this paper a \textbf{global tensor functor}
$\rkf_Q$ whose properties we summarize here:
\pagebreak
\begin{theoremnonum}
Let $V$ be a representation of a connected quiver $Q$. The functor 
\[
\rkf_Q: \repq \to \repq
\]
has the following properties:
\begin{enumerate}[(a)]
\item $\rkf_Q$ is a tensor functor, as defined above, and commutes with duality. (Theorem \ref{rankmultthm}, Proposition \ref{rankfirstprop})
\item For every arrow $a$ of $Q$, the linear map $(\rkf_Q V)_a$ is an isomorphism. Thus we can define a \textbf{global rank functor} on $Q$ by $\rankfunc_Q :=\imfuncv{x} \! \circ \, \rkf_Q$, whose isomorphism class is independent of $x$.  (Proposition \ref{rankfirstprop})
\item When $Q$ is a tree, $\rkf_Q V$ is isomorphic to a direct summand of $V$. (Theorem \ref{treesummandthm})
\end{enumerate}
\end{theoremnonum}

Via restriction, the global rank functors of subquivers $P \subseteq Q$ give rank functors on $Q$; for example, when $P$ is just a path $p$, considered as a subquiver of $Q$, we get the above functors $\imfunc_p$.  More generally, for any map of directed graphs $\alpha : Q' \to Q$,  the pushforward of the global rank functor of $Q'$ along $\alpha$ is a rank functor on $Q$.
Sometimes, a well-chosen $\alpha$ gives a rank functor on $Q$ that does not come from the global rank functor of any subquiver of $Q$ (Section \ref{ranktensorsect}).  In this paper, we focus on two running examples to illustrate how functions on the representation ring, obtained from rank functors, can be used to determine structure of $R(Q)$.

\subsection*{Acknowledgement}  The author wishes to thank his advisor Harm Derksen for many valuable discussions and support.

\section{Definitions}\label{defssect}
Throughout, $Q = (Q_0, Q_1, t, h)$ is a quiver on a finite vertex set $Q_0$ with finite arrow set $Q_1$.  The maps 
\[
t, h : Q_1 \to Q_0
\]
give the ``tail'' and ``head'' of an arrow, respectively.  We allow $Q$ to have oriented cycles, but  for simplicity we will assume that $Q$ is connected.  A \textbf{subquiver} $P \subseteq Q$ will also always be assumed to be connected.  We fix a field $K$ of any characteristic.
A \textbf{representation} of a quiver $Q$ is a collection of finite dimensional $K$-vector spaces
$\{V_x\}_{x \in Q_0}$,
and linear maps
$\{V_a: V_{ta} \to V_{ha}\}_{ a \in Q_1}$.
When $p= a_n \cdots a_2 a_1$ is a \textbf{path} in $Q$, we write $V_p:= V_{a_n} \cdots V_{a_2} V_{a_1}$.
A \textbf{morphism} $\varphi : V \to W$ between representations of a quiver $Q$ is given by specifying a linear map at each vertex
\[
\{ \varphi_x : V_x \to W_x \}_{x \in Q_0}
\]
such that these maps commute with the maps assigned to the arrows in $V$ and $W$, that is,
\[
\varphi_{ha} \circ V_{a} = W_{a} \circ \varphi_{ta}
\]
for $a \in Q_1$. We denote by $\repq$ the category of representations of $Q$.
The \textbf{dimension vector} of a representation $V$, written $\dimv V \in \N^{Q_0}$, is defined by $(\dimv V) _{x} := \dim_K V_{x}$, and the \textbf{support} of $V$ is the set 
\[
\supp V := \setst{x \in Q_0}{V_{x} \neq 0} .
\]
We say that $Q$ is a \textbf{tree} if the underlying graph is a tree (i.e., if removing any edge makes the graph disconnected).
The \textbf{opposite quiver} $Q^{op}$ of a quiver $Q = (Q_0, Q_1, t, h)$ is given by reversing the orientation of all arrows, so $Q^{op} = (Q_0, Q_1, h, t)$.  The categories $\repq^{op}$ and $\rep(Q^{op})$ are equivalent.
Vector space duals are denoted by superscript $^*$, and the \textbf{duality functor} on $\repq$ by 
\[
\begin{split}
D: \repq \to \rep (Q^{op}) \\
\left(DV\right) _{x} = V_{x}^*  \qquad \left(DV\right) _{a} = V_{a}^*
\end{split}
\]
for $V \in \repq$, $x \in Q_0$ and $a \in Q_1$.

\subsection{Tensor products of quiver representations}\label{tensorsect}

\begin{definition}
The \textbf{tensor product} of two quiver representations $V, W \in \repq$ is defined ``pointwise'':
\begin{align*}
(V \otimes W)_{x} &:= V_{x} \otimes W_{x} \qquad  x \in Q_0 \\
(V \otimes W)_{a} &:= V_{a} \otimes W_{a} \qquad a \in Q_1 .
\end{align*}
\end{definition}

An interesting, more general discussion of the tensor product of quiver representations can be found in \cite{herschtensorjpaa}.  We introduce the following notation.
\begin{definition}
For any quiver $Q$, we define the \textbf{identity representation} $\id_Q$ of $Q$ by
\[
(\id_Q)_{x} = K \qquad \qquad (\id_Q)_{a} = id_K
\]
for all $x \in Q_0$ and $a \in Q_1$.  (The subscript $Q$ is often omitted.)
\end{definition}

Note that $\id \otimes V \cong V$ for any representation $V$.
For the reader's convenience, we list some similarities and differences between the tensor product of quiver representations and tensor product of vector spaces.  A more complete treatment, in the language of tensor categories, will be given in Section \ref{categorysect}.

\begin{enumerate}[(a)]
\item Quiver tensor product is is a bifunctor
\[
\otimes : \repq \times \repq \to \repq
\]
which is associative, commutative, and distributes over direct sum of quiver representations.
\comment{
\[
V \otimes W \cong W \otimes V \qquad \  \text{and} \  \qquad  V \otimes (W \oplus U) \cong (V \otimes W) \oplus (V \otimes U)
\]
}
\item If we fix a representation $V$, the functor $T_V : \repq \to \repq$ defined by setting $T_V (W) = V \otimes W$ is exact.  This follows from the fact that exactness of a sequence of morphisms of quiver representations can be checked at each vertex. However, $T_V$ is not faithful, in general. For example, when $Q$ has 
no oriented cycles, the tensor product of two non-isomorphic simple representations is 0.

\item Quiver tensor product commutes with duality: $D (V \otimes W) \cong DV \otimes DW$.

\end{enumerate}

Although we know from Kac's Theorem \cite{kacrootsystemsa} that the dimension vectors of the indecomposable representations of a quiver $Q$ do not depend on the orientation of $Q$, the tensor product structure of $\repq$ \emph{does} depend on the orientation of $Q$.

\begin{notation}
We often denote a representation of a quiver by its dimension vector, if that representation is in the unique indecomposable isomorphism class of that dimension.
\end{notation}

\begin{example}\label{tensororienteg1}
Let $Q$ be the 3-subspace quiver,
\[
Q = \threesubspaceq
\]
and let $V$ be the indecomposable representation
\[
V = \threesubspacemaps{K^2}{K}{K}{K}{A}{B}{C} \qquad A={\eone}\quad B={\etwo} \quad C={\twobyone{1}{1}} .
\]
Then we can calculate
\[
V \otimes V \simeq \threesubspaced{1}{1}{0}{0} \oplus \threesubspaced{1}{0}{1}{0} \oplus \threesubspaced{1}{0}{0}{1} \oplus \threesubspaced{1}{0}{0}{0}
\]
since $\{m, n \}$ linearly independent in some vector space $L$ implies that $\{m \otimes m,\, n \otimes n,\, (m+n)\otimes (m+n) \}$ is linearly independent in $L \otimes L$.
\end{example}

\begin{example}\label{tensororienteg2}
However, if we change the orientation by flipping one of the arrows,
\[
Q = \QAq
\]
and let $W$ be the indecomposable representation of the same dimension vector as before,
\[
W = \QAmaps{K}{K}{K^2}{K}{A}{B}{C} \qquad A={\eone}\quad B={\etwo} \quad C={\onebytwo{1}{1}}
\]
then in this case we find that
\[
W \otimes W \simeq W \oplus \QAd{0}{0}{1}{0} \oplus \QAd{0}{0}{1}{0} .
\]
This can be directly calculated by writing down matrices, taking their tensor products, then finding the correct change of basis to put them in block form.  However, in Section \ref{ranktensorsect}, we will use rank functors to determine this decomposition without calculating any tensor products or change of basis.
\end{example}

\subsection{Tensor categories}\label{categorysect}
The purpose of this paper is to construct certain functors that are useful in studying tensor products of quiver representations.  The relevant properties of these functors can be stated concisely using the language of tensor categories, which we summarize here.
Although this subsection is necessary for technical purposes, we hope that it will not obscure the main ideas of the paper, which can be understood purely in terms of quiver representations.


The category $\repq$ is an abelian $K$-category \cite[p.~407--409]{assemetal}, so the spaces $\Hom_Q (V, W)$ are $K$-vector spaces,  and we will be interested in functors that preserve this structure.
We will use the following characterization: a functor
\[
F: \repq \to \rep(Q')
\]
is \textbf{additive} if and only if it \textbf{preserves direct sums}, meaning that each direct sum in $\repq$ with  insertion maps $i_V, i_W$ and projection maps $p_V, p_W$,
\[
\binarybiproduct{V}{V \oplus W}{W}{i_V}{i_W}{p_V}{p_W},
\]
is taken to an isomorphism in $\rep(Q')$
\[
F(V) \oplus F(W) \cong F(V \oplus W)
\]
with insertion maps $F(i_V),\, F(i_W)$ and projection maps $F(p_V),\, F(p_W)$
\cite[p.~67]{popescuabeliancats}.
Subfunctors and quotient functors of additive functors are additive \cite[p.~81]{popescuabeliancats}.

The additive bifunctor $\otimes$ (along with identity object $\id$) endows $\repq$ with the structure of a (relaxed) symmetric monoidal category \cite[VII]{Mcategories}, or simply a \textbf{tensor category} \cite{DMtannakian}.  This amounts to saying that the tensor product is functorial, satisfying some associativity and commutativity conditions, and has an identity object.
We summarize the definition of a tensor category, following \cite[p.~104--105]{DMtannakian}, but omitting some technicalities relating to associativity.  Consider a pair $(\categoryc, \otimes)$, where \categoryc\ is a category and $\otimes$ is a functor
\[
\otimes : \categoryc \times \categoryc \to \categoryc \qquad (X, Y) \mapsto X \otimes Y .
\]
Now let $\phi$ and $\psi$ be functorial isomorphisms
\[
\phi_{X, Y, Z} : X \otimes ( Y \otimes Z) \xto{\cong} (X \otimes Y ) \otimes Z \qquad \psi_{X, Y} : X \otimes Y \xto{\cong} Y \otimes X .
\]
We say that $\phi$ is an \textbf{associativity constraint} for $(\categoryc, \otimes)$ if $\phi$ satisfies a ``pentagon axiom'', and that $\psi$ is a \textbf{commutativity constraint} if $\psi_{Y, X} \circ \psi_{X, Y} = id_{X \otimes Y}$.  Such constraints are compatible with one another if they satisfy a ``hexagon axiom''.  The axioms omitted here are that certain diagrams of functorial isomorphisms involving $\phi$ and $\psi$ are commutative.
An \textbf{identity object} for $(\categoryc, \otimes)$ is an object $U$ of $\categoryc$ and an isomorphism
\[
 u : U \to U \otimes U
\]
such that the functor
\[
T_U : \categoryc \to \categoryc \qquad X \mapsto U \otimes X
\]
is an equivalence of categories.
A system $(\categoryc, \otimes, \phi, \psi)$ as above will be called a \textbf{tensor category} if the constraints are compatible and there is an identity object.

The category of finite dimensional vector spaces, with the standard tensor product and standard associativity and commutativity isomorphisms, is a tensor category.  For an arbitrary quiver $Q$, let us note once and for all that a morphism $\alpha = \{ \alpha_x \}_{x \in Q_0}$ is an isomorphism if and only if each $\alpha_x$ is an isomorphism, and similarly, commutativity of diagrams can be checked at each vertex of $Q$.  Then since our tensor product $\otimes$ is just the standard tensor product of finite dimensional vector spaces at each vertex, defining associativity and commutativity constraints pointwise equips $\repq$ with the structure of a tensor category (with identity object $\id$).

A \textbf{tensor functor} \cite[p.~113--114]{DMtannakian} is a pair $(F, c)$ consisting of an additive functor 
\[
F : \repq \to \rep (Q')
\]
and a functorial isomorphism
\[
c_{X, Y} : F(X) \otimes F(Y) \cong F(X \otimes Y)
\]
such that:
\begin{enumerate}[(a)]
\item The isomorphism $c$ is compatible with associativity, expressed by another ``hexagon axiom''.
\item $c$ is compatible with commutativity, that is, $c_{Y, X} \circ \psi_{F(X), F(Y)} \cong F(\psi_{X, Y}) \circ c_{X,Y}$.
\item $F(\id_Q) \cong \id_{Q'}$ .
\end{enumerate}

If we relax the condition that $c$ be an isomorphism by simply requiring the existence of either
\[
c_{X, Y} : F(X) \otimes F(Y) \to F(X \otimes Y) \qquad \text{or} \qquad c_{X, Y} : F(X \otimes Y) \to F(X) \otimes F(Y)
\]
satisfying (a), (b), and (c) (appropriately modified),
we will say that $(F, c)$ is a \textbf{weak tensor functor}.
Note that the definition of a tensor functor covers both covariant and contravariant functors, since all the maps appearing in the definition are isomorphisms.
For example, the duality functor $D$ is a tensor functor.

\begin{remark}\label{iteratedtensor}
The technical axioms involving the associativity and commutativity constraints give us the following facts about tensor categories and tensor functors \cite[VII.2]{Mcategories}, \cite[p.~106--108,~113--115]{DMtannakian}:
\begin{enumerate}[(a)]
\item There is an essentially unique way to extend the tensor product to any finite family of objects of \categoryc, 
\[
\tensor_{i =1}^k : \categoryc^k \to \categoryc \qquad (X_i ) \mapsto \tensor_{i =1}^k X_i .
\]
\item There is an action of the symmetric group $\symgp{k}$ on tensor products
\[
\sigma \cdot \tensor_{i =1}^k X_i := \tensor_{i =1}^k X_{\sigma (i)} \qquad \sigma \in \symgp{k} .
\]
\item For any tensor functor $F$, there is an isomorphism of functors
\[
\tensor_{i =1}^k F(X_i) \cong F \left( \tensor_{i =1}^k X_i \right)
\]
which is $\symgp{k}$-equivariant, that is,
\[
F \left( \sigma \cdot \tensor_{i =1}^k X_i \right) \cong \sigma \cdot \tensor_{i =1}^k F \left( X_{i} \right) .
\]
\end{enumerate}
\end{remark}

\begin{definition}
A \textbf{rank functor} on $Q$ is a tensor functor from $\repq$ to $\rep(A_1) = \vscat$.
\end{definition}

\subsection{Schur functors on quiver representations}
The Schur functors $\schurl$ \cite[p.~76]{fultonreptheory} act on quiver representations pointwise:
\begin{align*}
\left( \schurl V \right) _{x} &:= \schurl \left(V_{x} \right) \qquad x \in Q_0 \\
\left( \schurl V \right) _{a} &:= \schurl \left(V_{a} \right) \qquad a \in Q_1
\end{align*}
where $\lambda$ is a partition and $V \in \repq$. From the functoriality of $\schurl$ on vector spaces, it follows that this defines a functor from $\repq$ to $\repq$.

\begin{example}
Using $V$ from Example \ref{tensororienteg1}, we get
\[
S^2 V \simeq \threesubspaced{1}{1}{0}{0} \oplus \threesubspaced{1}{0}{1}{0} \oplus \threesubspaced{1}{0}{0}{1} \qquad \text{and } \qquad
S^{(1,1)} V = \exterior^2 V \simeq \threesubspaced{1}{0}{0}{0} .
\]
In particular, note that $S^2 V$ is not indecomposable.
\end{example}

Now suppose that $\charac K = 0$ in the remainder of this subsection.  For a vector space $V_x$, we have the Schur decomposition of a tensor power of $V_x$ \cite[p.~87]{fultonreptheory}:
\begin{equation}\label{schurdecompvseq}
\tensor^k V_x \simeq \bigoplus_{\lpartk} \schurl{(V_x)} \otimes \glamb
\end{equation}
where $\lpartk$ means that $\lambda$ is a partition of $k$, and $\glamb$ is the irreducible representation of the symmetric group $\symgp{k}$ corresponding to $\lambda$.  Because this is functorial in $V_x$, we expect to be able to utilize this decomposition for quiver representations.

\begin{definition}
A \textbf{$\qsymk$-representation} is a representation $V \in \repq$ such that each $V_{x}$ is a representation of $\symgp{k}$ and each map $V_{a}$ is $\symgp{k}$-equivariant.
A morphism of $\qsymk$-representations is just a morphism of quiver representations $\{ \varphi _{x} \}$ such that each $\varphi _{x}$ is $\symgp{k}$-equivariant.
\end{definition}

The $\qsymk$-representations form a category, so we have notions of subrepresentations, irreducible objects, and so forth.  

\begin{example}\label{qrepsymeg}
Any linear representation $V$ of $\symgp{k}$ gives rise to a $\qsymk$-representation $V^Q$ by setting $V^Q _{x} = V$ for all vertices $x$, and $V^Q_{a} = id$ for all arrows $a$.  Note that if $V$ is an irreducible $\symgp{k}$ representation, $V^Q$ is an irreducible $\qsymk$-representation, but $V^Q$ is not even necessarily indecomposable in $\repq$ or as a representation of $\symgp{k}$.  In fact, we have decompositions
\begin{equation}\label{VQdecompeq}
V^Q \simeq \bigoplus_{i=1}^{\dim V} \id \qquad \text{and} \qquad V^Q \simeq \bigoplus_{i=1}^{\# Q_0} V
\end{equation}
as a representation of $Q$, and as a representation of $\symgp{k}$, respectively.
\end{example}
 
If $V$ is any representation of $Q$, then $\tensor^k V$ becomes a $\qsymk$-representation by letting $\symgp{k}$ permute the factors.

\begin{prop}\label{schurdecompprop}
Let $Q$ be a quiver and $V \in \repq$.  Denote by $\glamb$ the irreducible linear representation of $\symgp{k}$ corresponding to $\lambda$.  Then we have a direct sum decomposition
\begin{equation}\label{schurdecompqeq}
\tensor^k V \simeq \bigoplus_{\lpartk} \schurl{V} \otimes \glamb^Q
\end{equation}
as $\qsymk$-representations, which is functorial in $V$.
\end{prop}
\begin{proof}
The isomorphism is defined at each vertex by (\ref{schurdecompvseq}).  Functoriality of (\ref{schurdecompvseq}) in $V_x$ implies that (\ref{schurdecompqeq}) is an isomorphism in $\repq$.  Since $\symgp{k}$ acts trivially on each $\schurl{V}_{x}$, and the identity map is $\symgp{k}$-equivariant, the induced maps $\bigoplus_{\lpartk} \schurl{V_{a}} \otimes id_{\glamb}$ are evidently $\symgp{k}$-equivariant, so the right hand side is in fact a $\qsymk$-representation.  Then because (\ref{schurdecompvseq}) is an isomorphism of $\symgp{k}$ representations, (\ref{schurdecompqeq}) is an isomorphism of $\qsymk$-representations.  Functoriality in $V$ follows from the functoriality at each vertex of (\ref{schurdecompvseq}).
\end{proof}

\begin{theorem}\label{rankschurthm}
Let $\charac K = 0$. Then any tensor functor $F : \repq \to \repq$ commutes with the Schur functors.  That is, there is an isomorphism of functors
\begin{equation*}
F \circ \schurm{} \cong \schurm{} \circ F .
\end{equation*}
\end{theorem}

\begin{proof}
We demonstrate the isomorphism for an object $V \in \repq$; functoriality in $V$ will be clear at each step of the proof, without explicit mention.  Let $k = |\mu|$.
By Remark \ref{iteratedtensor}, there is an isomorphism of $\qsymk$-representations
\begin{equation}\label{rankmultkeq}
F (\tensor^k V ) \cong \tensor^k F V .
\end{equation}
To prove the theorem, we just write down the Schur decomposition of each side, then try to match up the correct pieces.

Applying Proposition \ref{schurdecompprop} to the left hand side of (\ref{rankmultkeq}), we have isomorphisms of $\qsymk$-representations
\begin{equation*}
F \left( \tensor^k V \right) \cong F \left( \bigoplus_{\lpartk} \schurl{V} \otimes \glamb^Q \right) \cong \bigoplus_{\lpartk} F \left( \schurl{V} \right) \otimes F(\glamb^Q)  \cong \bigoplus_{\lpartk} F \left( \schurl{V} \right) \otimes \glamb^Q .
\end{equation*}
Note that $F (\glamb^Q) = \glamb^Q$ by the first isomorphism of (\ref{VQdecompeq}).

Then applying Proposition \ref{schurdecompprop} to the right hand side of (\ref{rankmultkeq}), we have an isomorphism of $\qsymk$-representations
\begin{equation*}
\tensor^k F V \cong \bigoplus_{\lpartk} \schurl{\left( F V \right)} \otimes \glamb^Q
\end{equation*}
hence an isomorphism of $\qsymk$-representations
\begin{equation*}
\bigoplus_{\lpartk} F \left( \schurl{V} \right) \otimes \glamb^Q \cong \bigoplus_{\lpartk} \schurl{\left( F V \right)} \otimes \glamb^Q .
\end{equation*}

Each $\glamb$-isotypic component of the left hand side must map to the $\glamb$-isotypic component of the right hand side, so by dimension count
we get isomorphisms of the summands on each side indexed by the same partition $\lambda$.  This gives an isomorphism in $\repq$
\begin{equation*}
\bigoplus_{\dim G_\mu} F \left( \schurm{V} \right) \simeq F \left( \schurm{V} \right) \otimes G_\mu^Q \cong \schurm{\left( F V \right)} \otimes G_\mu^Q \simeq \bigoplus_{\dim G_\mu} \schurm{\left( F V \right)}
\end{equation*}
and so by the Krull-Schmidt property of $\repq$, we have
\[
F \left( \schurm{V} \right) \cong \schurm{\left( F V \right)} . \qedhere
\]
\end{proof}

\section{The Representation Ring of a Quiver}\label{repringsect}

To study tensor products of, say, finite dimensional complex representations of $SL_n$, one starts by describing the indecomposable representations of $SL_n$.  In this case, they are indexed by partitions $\lambda$, and one can study tensor products of representations of $SL_n$ via combinatorics of these partitions \cite{fultonreptheory}.

For an arbitrary quiver $Q$, however, there is no good description of the indecomposable representations of $Q$ to use as a starting point.  Alternatively, we can look at the representation ring of $Q$ and study the abstract properties of this ring.  This gives us a convenient setting for formulating properties of tensor products.

\begin{definition}
Let $[V]$ denote the isomorphism class of a representation $V$.  Then define $R(Q)$ to be the free abelian group generated by isomorphism classes of representations of $Q$, modulo the subgroup generated by all $[V \oplus W] - [V] - [W]$.  The operation
\[
[V] \cdot [W] := [V \otimes W] \qquad \text{for }V,\ W \in \repq
\]
induces a well-defined multiplication on $R(Q)$, making $R(Q)$ into a commutative ring with identity $[\id_Q]$, called the \textbf{representation ring} of $Q$.  We usually omit the brackets $[\ ]$ and just refer to representations of $Q$ as elements of $R(Q)$.
\end{definition}

Although we introduce virtual representations to form this ring, every element $r \in R(Q)$ can be written as a formal difference
\begin{equation*}
r = V - W \qquad \text{with }V,\ W \in \repq .
\end{equation*}
Then any additive (resp. multiplicative) relation $z = x + y$ (resp. $z = xy$) can be rewritten to give some isomorphism of actual representations of $Q$.  

The Grothendieck group of $\repq$ \cite[p.~87]{assemetal} is $R(Q) / \mathfrak{a}$, where $\mathfrak{a}$ is the subgroup generated by elements
\[
[V] - [U] - [W]
\]
for all short exact sequences
\[
\ses{U}{V}{W} . 
\]
We do not work modulo short exact sequences, because this loses too much information in our setting. For example, the Grothendieck group of a quiver $Q$ without oriented cycles is always isomorphic to $\Z^{Q_0}$, and the image of a representation in the Grothendieck group is just its dimension vector.

\begin{remark}
The ring $R(Q)$ generally depends on the base field $K$.
\end{remark}

\begin{remark}
If $\relideal$ is an ideal of relations in $Q$ \cite[II.2]{assemetal}, and $V, W$ are representations of the bound quiver $(Q, \relideal)$, then $V \otimes W$ may \emph{not} be a representation of $(Q, \relideal)$. 
For example, let $(Q, \relideal)$ be the bound quiver
\[
Q= \kronthreemaps{\bullet}{\bullet}{a}{b}{c} \qquad \relideal = \gen{a+b -c}
\]
and assume $\charac K \neq 2$.  Then for  the representation
\[
V = \kronthreemaps{K}{K}{(1)}{(1)}{(2)}
\]
we have that $V_c^{\otimes 2} = (4) \neq V_a^{\otimes 2} + V_b^{\otimes 2} = (1) + (1) = (2)$, so $V \otimes V \notin \rep(Q, \relideal)$.


However, if $\relideal$ is generated by commutativity relations (that is, relations of the form $p-q$ for paths $p, q$) then the representations of $(Q, \relideal)$ do generate a subalgebra of $R(Q)$.  If $\relideal$ is generated by zero relations (relations of the form $p = 0$ for $p$ a path), then representations of $(Q, \relideal)$ do not generate a subalgebra because $\id \notin \rep (Q, \relideal)$.  But because the tensor product of any map with a zero map is zero, these representations generate an ideal in $R(Q)$.
\end{remark}

We note some properties of the representation ring.
\begin{enumerate}[(a)]
\item By the Krull-Schmidt theorem \cite[Theorem~I.4.10]{assemetal}, $R(Q)$ is a free $\Z$-module, with basis given by the isomorphism classes of indecomposable representations of $Q$.
\item By Gabriel's theorem \cite{gabriel}, $R(Q)$ is a finitely generated $\Z$-module if and only if $Q$ is a Dynkin quiver.  In this case, 
\begin{equation*}
\rank_\Z R(Q) = \#\{\text{indecomposable representations of }Q \} = \#\{\text{positive roots of }Q \} .
\end{equation*}
\item Write $[\repq] := \setst{ [V] \in R(Q)}{V \in \repq}$ for the semiring of ``actual representations'' in $R(Q)$. The ring $R(Q)$ has the following universal property: if $A$ is a ring, and $f: [\repq] \to A$ a map of semirings, then $f$ uniquely extends to a ring homomorphism $f: R(Q) \to A$ (which we usually denote by the same symbol). That is, the functor sending the semiring $[\repq]$ to the ring $R(Q)$ is left adjoint to the forgetful functor from rings to semirings.
\end{enumerate}

\begin{remark}
The notion of $\lambda$-rings was introduced by Grothendieck in \cite{Gchernclasses}.  The idea is to define unary operations $\lambda^i$ on a commutative ring with identity, called ``$\lambda$-operations'', which have the formal properties of exterior power operations.  Essentially, one wants to express $\lambda^i (x+y),\, \lambda^i (xy)$, and $\lambda^i (\lambda^j (x))$ as some universal polynomials in the values of $\lambda^k (x)$ and $\lambda^k (y)$.  The reader is referred to \cite{knutsonlambda} and \cite{atiyahtall} for the definitions and basic properties.
 For example, the ring of integers is a $\lambda$-ring under the operations
\[
\lambda^i : \Z \to \Z \qquad n \mapsto \binom{n}{i} .
\]
  
In the quiver setting, we can define $\lambda$-operations on $R(Q)$ by setting
\[
\lambda^i V = \exterior^i V
\]
for $V \in \repq$.  For example, this gives the same $\lambda$-ring structure as above on $R(A_1) \cong \Z$, since $\dim \exterior^i V = \binom{\dim V}{i}$.  In the general case, because the exterior power operations on a quiver representation act as the standard exterior powers at each vertex, it is immediate that these operations give $R(Q)$ the structure of a $\lambda$-ring.

A homomorphism of $\lambda$-rings is just a ring homomorphism that commutes with the $\lambda$-operations.  In Theorem \ref{rankschurthm} we saw that when $\charac K =0$, tensor functors commute with the Schur operations on quiver representations, so in particular they commute with exterior powers.  Hence a homomorphism of representation rings induced from a tensor functor (in characteristic 0) is a $\lambda$-ring homomorphism.
\end{remark}

\section{Construction of the Global Tensor Functor}\label{constructionsect}

We can construct a canonical tensor functor $\rkf_Q : \repq \to \repq$ for any quiver $Q$.

\begin{definition}
Call a representation $V$ of a quiver $Q$ \textbf{epimorphic} (resp. \textbf{monomorphic}) if $V_{a}$ is surjective (resp. injective) for each arrow $a \in Q_1$.
\end{definition}

\begin{example}
If $Q$ has no oriented cycles, then injective representations are epimorphic and projective representations are monomorphic. Let $P$ be the projective representation corresponding to a vertex $x$, and $a \in Q_1$. For each vertex $y$, the vector space $P_y$ has a basis given by all paths from $x$ to $y$, and the maps $P_{a} : P_{ta} \to P_{ha}$ are given by composition with $a$ (cf. \cite[p.~79]{assemetal}).  If $p,\, q$ are distinct paths from $x$ to $ta$, then $ap,\, aq$ are distinct paths from $x$ to $ha$, so $P_a$ takes the standard basis of $P_{ta}$ to a subset of the standard basis of $P_{ha}$. Hence each map $P_{a}$ is injective.  The case of injective representations is similar.
\end{example}

For $V \in \repq$, the sum of any collection of epimorphic subrepresentations of $V$ is epimorphic, hence $V$ has a unique maximal epimorphic subrepresentation $\surjrep_Q (V)$.  
Dually, $V$ also has a maximal monomorphic quotient $\injrep_Q (V)$, and these are related by a canonical isomorphism of $Q^{op}$ representations
\begin{equation}\label{surjinjdualeq}
D \injrep_Q (V) \cong \surjrep_{Q^{op}} (DV ) .
\end{equation}

\begin{example}\label{surjinjeg}
Let $Q$ and $V$ be as in Example \ref{tensororienteg2}.
Then we get $\surjrep_Q (V) = 0$, essentially because $\im A \cap \im B = 0$,
and $\injrep_Q (V) \simeq \QAd{1}{1}{1}{1}$, given by $K^2 / \ker C$ at the branch point.
\end{example}

\begin{example}
Let $Q$ be the single loop quiver, so $V \in \repq$ is given by a vector space $V_0$ together with an endomorphism $A$.  Then if $V$ is \emph{indecomposable},
\[
\surjrep_Q (V) = \injrep_Q (V) = \begin{cases}
V & \text{if } A \text{ is an isomorphism} \\
0 & \text{otherwise}
\end{cases} .
\]
\end{example}

It is easy to check that for any quiver $Q$, both $\surjrep_Q$ and $\injrep_Q$ are covariant functors from $\repq$ to $\repq$, and (\ref{surjinjdualeq}) is a natural isomorphism of functors.  Furthermore, $\surjrep_Q$ (resp. $\injrep_Q$) is a sub- (resp. quotient-) functor of the identity functor, hence each is additive, and there is a natural transformation given by the composition
\[
\Phi : \surjrep_Q \into id_{\repq} \onto \injrep_Q .
\]

\begin{definition}
The \textbf{global tensor functor} of $Q$ is defined to be the image functor of this natural transformation:
\[
\rkf_Q := \im (\Phi) : \repq \to \repq .
\]
(The subscript $Q$ and parentheses around the input are often omitted from all of the above functors.)
\end{definition}

We can immediately note some useful properties of $\rkf_Q$.

\begin{prop}\label{rankfirstprop}  The global tensor functor is additive and commutes with duality. So we have natural isomorphisms of functors
\[
\rkf_Q V \oplus \rkf_Q W \cong \rkf_Q (V \oplus W) \qquad \text{and} \qquad D \circ \rkf_Q \cong \rkf_{Q^{op}} \circ D .
\]
Furthermore, for any $V \in \repq$, the representation $\rkf_Q V$ is both epimorphic and monomorphic, i.e., the linear maps $(\rkf_Q V)_{a}$ are isomorphisms for all $a \in Q_1$.
\end{prop}
\begin{proof}
The functor $\rkf$ is additive because it is a quotient of the additive functor $\surjrep$.
That $\rkf$ commutes with duality follows from (\ref{surjinjdualeq}) and the universal property of an image.
For the last statement, we have maps
\[
\surjrep V \onto \rkf V \into \injrep V .
\]
The conditions to be a morphism of quiver representations imply that a quotient of an epimorphic representation is epimorphic, and a subrepresentation of a monomorphic representation is monomorphic, hence $\rkf V$ is both.
\end{proof}

In particular, the dimension of $(\rkf_Q V)_{x}$ is independent of $x \in Q_0$ when $Q$ is connected (as all quivers in this paper are). We will later prove (Theorem \ref{rankmultthm}) that $\rkf_Q$ is actually a tensor functor, which leads us to define:

\begin{definition}
The \textbf{global rank functor} of a quiver $Q$ is
\[
\rankfunc_Q := \imfuncv{x} \! \circ \, \rkf_Q : \repq \to \vscat.
\]
If $F$ is any rank functor on a quiver $Q$, then the induced ring homomorphism 
\[
f: R(Q) \to \Z
\]
is called the \textbf{rank function} associated to $F$.  Thus the \textbf{global rank function} of a connected quiver $Q$ is given by $r_Q (V) = \dim_K (\rkf_Q V)_{x}$ for $V \in \repq$, and extended by linearity to $R(Q)$.  By the remark above, this is independent of the choice of $x \in Q_0$.
\end{definition}

We can compute some simple examples of the global tensor functor.

\begin{example}
Let $Q$ be equioriented of type $A_3$,
\[
Q = \equiAthree{\bullet}{\bullet}{\bullet}{a}{b} .
\]
Then for a representation $V = \equiAthree{V_1}{V_2}{V_3}{A}{B}$, one can compute from the definitions that
\[
\rkf_Q V = \equiAthree{\frac{V_1}{\ker BA}}{\frac{\im A}{\ker B \cap \im A}}{\im BA}{\overline{A}}{\overline{B}}
\]
so the global rank functor is $\rankfunc_Q \cong \imfunc_{ba}$, and the associated rank function is $r_Q (V) = \rank (BA)$.
This easily generalizes to a quiver of type equioriented $A_n$.
\end{example}

The global rank function does not always correspond to the rank of some map, but can still be easily described sometimes.

\begin{example}
When $Q$ is the two subspace quiver
\[
Q: \twosubspaceq
\]
we can again explicitly compute the global tensor functor.  If $V = \twosubspacemaps{V_1}{V_2}{V_3}{A}{B}$, then
\[
\rkf_Q V = \twosubspacemaps{\frac{A^{-1}(\im B)}{\ker A}}{ \im A \cap \im B}{\frac{B^{-1}(\im A)}{\ker B}}{\overline{A}}{\overline{B}}
\]
and so $r_Q (V)= \dim_K (\im A \cap \im B)$.
\end{example}

When $Q$ has many sinks and sources, the global rank function becomes more cumbersome to write down explicitly:
\begin{example}
Let $Q$ be of type $A_4$ and the alternating orientation
\[
Q = \alterAfourq .
\]
If we write a representation as 
\[
V =\alterAfour{V_1}{V_2}{V_3}{V_4}{A}{B}{C}
\]
then we can compute from the definitions that
\[
r_Q (V) = \dim_K \left( \frac{\im A \cap \im  B}{ \im A \cap B (\ker C) } \right) = \dim_K \left( \frac{ B^{-1}(\im A ) }{ \ker B + \ker C } \right) . 
\]
\end{example}

\begin{remark}
Since $\rkf V$ is both a monomorphic quotient of $\surjrep V$, and an epimorphic subrepresentation of $\injrep V$, the universal properties yield natural transformations
\[
\injrep \circ \surjrep \onto \rkf \into \surjrep \circ \injrep ,
\]
Neither of these are necessarily isomorphisms: let $Q$ and $V$ be as in Example \ref{surjinjeg}.  Then 
\[
\rkf V = 0 \qquad \text{but} \quad \surjrep \left( \injrep V \right) \simeq \QAd{1}{1}{1}{1} .
\]
Dualize to get an analogous example for the other map. (We will not be interested in compositions of $\surjrep, \rkf,$ and $\injrep$ with one another in this paper.)

\end{remark}

\section{Global Tensor Functors for Trees}

A quiver $Q$ generates a category $\quivcat(Q)$ by taking the objects of $\quivcat(Q)$ to be the vertices of $Q$, and the morphisms of $\quivcat(Q)$ to be the paths in $Q$:
\[
\Ob \quivcat(Q) := Q_0 \qquad \Mor_{\quivcat(Q)} (x, y) := \{ \text{paths from $x$ to $y$} \} .
\]
The trivial path at a vertex $x$ is the identity morphism for $x$, and composition of morphisms is composition of paths.  A representation $V$ of $Q$ is the same thing as a functor from $\quivcat(Q)$ to $\vscat$, denoted $\repfunc_V$, and a morphism of representations is a natural transformation of the corresponding functors.  In other words, a quiver representation is a \textbf{diagram} of type $\quivcat(Q)$ in $\vscat$.

Taking the limit and colimit of such a diagram $\repfunc_V$, we get vector spaces $\plim V$ and $\ilim V$, respectively, with natural maps
\[
\alpha_x : \plim V \to V_{x} \qquad \text{and} \qquad \beta_x : V_{x} \to \ilim V
\]
for each $x \in Q_0$.  These maps satisfy $\alpha_{ha} = V_a \circ \alpha_{ta}$ and $\beta_{ta} = \beta_{ha} \circ V_a$ for every arrow $a \in Q_1$,
and therefore $\eta_V := \beta_x \circ \alpha_x$ does not depend on $x$.  When $Q$ is a tree, we will see that the functors $\surjrep, \injrep,$ and $\rkf$ can be constructed using limits, and have a nice connection to $\Hom$ spaces in $\repq$. 

\begin{prop}\label{limpairprop} 
There are functorial isomorphisms of vector spaces
\begin{equation*}
\plim V \cong \Hom_Q (\id , V) \qquad \text{and} \qquad \ilim V \cong \Hom_Q (V, \id )^* ,
\end{equation*}
so $\plim V$ is representable by $\id$.
The natural map $\eta_V : \plim V \to \ilim V$ corresponds to the pairing
\[
\Hom_Q (\id, V ) \times \Hom_Q (V, \id) \to \Hom_Q (\id, \id ) \cong K \qquad (f, g) \mapsto g \circ f .
\]
\end{prop}
\begin{proof}
Fix a compatible basis $\{e_x \}_{x \in Q_0}$ of $\id$, that is, vectors $e_x \in \id_{x} \simeq K$ such that 
\[
\id_{a}(e_{ta}) = e_{ha}
\]
for every arrow $a$.  Now given $v \in \plim V$, define $f_v \in \Hom_Q (\id , V)$ by $(f_v)_{x}(e_x) = \alpha_x (v)$.  It is easy to check from definitions and universal properties that $f_v$ is a $Q$-morphism, that $v \mapsto f_v$ gives a vector space isomorphism $\plim V \xto{\sim} \Hom_Q (\id, V)$, and that this isomorphism is natural in $V$.

By applying this to $Q^{op}$ and dualizing, we get $\ilim V \cong \Hom_Q (V, \id )^*$.  With this, it is routine to check that $\eta_V$ corresponds to the stated natural pairing.
\end{proof}

The following proposition relates these spaces to the global tensor functor when $Q$ is a tree.

\begin{prop}\label{surjinjlimitprop}
If $Q$ is a tree, we can construct the functors $\surjrep$ and $\injrep$ from limits and colimits:
\begin{equation*}
(\surjrep V)_x = \im \alpha_x \qquad \text{and} \qquad (\injrep V)_x = \frac{V_{x}}{\ker \beta_x} 
\end{equation*}
where $\alpha_x$ and $\beta_x$ are defined for $V \in \repq$ above, and the maps $(\surjrep V)_{a}$, $(\injrep V)_{a}$ are induced from $V_{a}$.
Thus, for each $x \in Q_0$, we have
\[
(\rkf V)_x = \frac{\im \alpha_x}{\ker \beta_x \cap \im \alpha_x} \cong \im \eta_V \qquad \text{and so} \quad r_Q (V) = \rank \eta_V .
\]
\end{prop}

\begin{proof}
For each arrow $a \in Q_1$, the universal property of $\plim$ gives a commutative diagram:
\[
\xymatrix@R=2ex{
	&	{\im} \alpha_{ta} \subseteq V_{ta} \ar@<4ex>[dd]^{V_{a}} \ar@{-->}@<-2ex>[dd] \\
{\plim V} \ar@{->>}[ur]^{\alpha_{ta}} \ar@{->>}[dr]_{\alpha_{ha}} \\
	&	{\im} \alpha_{ha} \subseteq V_{ha}
}
\]
which shows that $N := \bigoplus_{x \in Q_0} \im \alpha_x$ is an epimorphic subrepresentation of $V$.
We will show that any epimorphic subrepresentation of $V$ is contained in $N$.

Now let $E \subseteq V$ be an arbitrary epimorphic subrepresentation of $V$.  For any vertex $x$, and $v \in E_{x}$, we claim that there is some $f \in \Hom_Q (\id, E)$ such that $v \in \im f$.  This is proved by induction on the number of vertices of $Q$, using the notation of Proposition \ref{limpairprop}.  If $Q$ has one vertex, the situation trivial.  Otherwise, choose a vertex $y \neq x$ such that there is precisely one arrow $a \in Q_1$ with $y = ta$ or $y = ha$.  This is possible because a tree always has at least two such vertices.  Let $P \subset Q$ be the subquiver of $Q$ obtained by removing $y$ and $a$, so that $x \in P_0$ and $v \in \left( E |_P \right)_{x}$, where $E|_P$ denotes the restriction of $E$ to $P$.  Then by induction, there exists $f \in \Hom_P \left( \id_P, E|_P \right)$ such that $v \in \im f$.

We can extend $f$ to $Q$: if $y = ha$, then simply set $f(e_y) = E_{a}\left( f(e_{ta}) \right)$. If $y = ta$, then since $E$ is epimorphic, there exists some $w \in E(y)$ such that $V_{a}(w) = f(e_{ha})$.  In this case, set $f(e_y) = w$; it is immediate from the definition that in either situation $f \in \Hom_Q (\id_Q, E)$.

Regarding $f \in \Hom_Q (\id_Q, V)$ via the inclusion $E \subseteq V$, the explicit formulation of the isomorphism $\plim V \cong \Hom_Q (\id, V)$ in the proof of Proposition \ref{limpairprop} shows that $v \in \im \alpha_x$, hence $v \in N$.  Hence $E \subseteq N$.  So $N$ must be maximal epimorphic, that is, $N = \surjrep V$.

The equation for $\injrep$ follows by dualizing, then the equation for $\rkf$ from the other two equations and its definition.
Since these are equalities as a subrepresentation and quotient representation of $V$, respectively, we get the following commutative diagram:
\[
\xymatrix@=3ex{
{\plim V} \ar@{->>}[dr] \ar[rr]^{\alpha_x} \ar@/^4ex/[rrrr]^{\eta_V} & & V_{x} \ar@{->>}[dr] \ar[rr]^{\beta_x} & & {\ilim V} \\
	& {\im} \alpha_x = (\surjrep V)_{x} \ar@{^`->}[ur] \ar@{->>}[dr]& & {\im} \beta_x \cong (\injrep V)_{x} \ar@{^`->}[ur]\\
	&	& ({\rkf} V)_{x} \ar@{^`->}[ur]
}
\]
In particular note that $(\rkf V)_{x} \cong \im \eta_V$.
\end{proof}

This characterization of $\rkf_Q$ allows us to see that $\rkf_Q V$ is isomorphic to a direct summand of $V$ when $Q$ is a tree.

\begin{theorem}\label{treesummandthm}
Let $Q$ be a tree, and $V$ an indecomposable representation of $Q$.  Then
\[
\rkf V \neq 0 \iff V \simeq \id . 
\]
In particular, if we write
\[
V \simeq \bigoplus_j V_j
\]
where $V_j \subseteq V$ are indecomposable subrepresentations of $V$, then
\[
\rkf V \simeq \bigoplus_{V_j \simeq \id} V_j
\]
\end{theorem}
\begin{proof}
Suppose $V$ is an indecomposable representation of $Q$. If $\rkf V = 0$, then certainly $V \not \simeq \id$ because $\rkf \id = \id$.  If $\rkf V \neq 0$, then by Proposition \ref{surjinjlimitprop}, $\im \eta_V \neq 0$.  Then Proposition \ref{limpairprop} implies that there is a pair $(f, g) \in \Hom_Q (\id, V ) \times \Hom_Q (V, \id)$ such that
\[
\xymatrix{ {\id} \ar[r]_{f} \ar@/^2ex/[rr]^{id} &	V \ar[r]_{g}& {\id} 	}
\]
so $V$ has a direct summand isomorphic to $\id$.  But we took $V$ to be indecomposable, so $V \simeq \id$.  The second statement follows from additivity of $\rkf$ (Proposition \ref{rankfirstprop}).
\end{proof}

\begin{example}
Let $Q$ be the following quiver, and $V$ an indecomposable representation of $Q$.
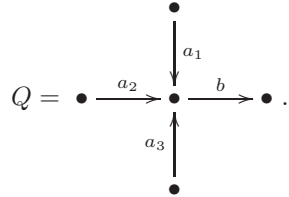
\[
Q = \extDfouramaps{\bullet}{\bullet}{\bullet}{\bullet}{\bullet}{a_1}{a_2}{a_3}{b} . 
\]
Also suppose that $V$ is not simple, so that each $V_{a_i}$ is injective and $V_b$ is surjective.  Denoting the branch vertex of $Q$ by $x$, one can calculate from the definitions that
\[
(\surjrep V)_{x} = \bigcap_{i = 1}^3 V_{a_i} \qquad \qquad (\injrep V)_{x} = \frac{V_{x}}{\ker V_b}
\]
so by Theorem \ref{treesummandthm} we have
\[
V \simeq \id \iff (\rkf_Q V)_{x} \neq 0 \iff \frac{\bigcap\limits_{i = 1}^3 V_{a_i} + \ker V_b}{\ker V_b} \neq 0 \iff  \bigcap_{i = 1}^3 V_{a_i} \nsubseteq \ker V_b . 
\]
\end{example}

One can check for some other quivers which are not trees, for example when $Q$ is the Kronecker quiver
\[
Q: \krontwo{\bullet}{\bullet} ,
\]
that $\rkf_Q V = 0$ when $V$ is indecomposable and some $V_{a}$ is not an isomorphism.  This is easy for this particular example because $Q$ is of tame representation type, so we have nice descriptions of the indecomposable representations of $Q$.  In this case we can still say that  $\rkf_Q$ ``picks out'' the indecomposable summands for which the map over every arrow is an isomorphism, although these representations are not necessarily isomorphic to $\id_Q$.
The following example, however, shows that this property does not hold for all quivers.

\begin{example}
This example shows that $\rkf_Q V$ is not necessarily isomorphic to a direct summand of $V$.
Let $Q$ be the generalized Kronecker quiver with four arrows $a, b, c, d$.
\[
Q: \kronfour{\bullet}{\bullet}
\]
Let $V$ be the representation with dimension vector $\alpha = (2, 3)$ and maps given by
\[
	V_{a} = \begin{pmatrix} 1& 1 \\ 0& 0 \\ 0& 0 \end{pmatrix}	
\qquad V_b = \begin{pmatrix} 1& 0 \\ 0& 1 \\ 0& 0 \end{pmatrix}	
\qquad V_c = \begin{pmatrix} 1& 0 \\ 0& 0 \\ 0& 1 \end{pmatrix}	
\qquad V_d = \begin{pmatrix} 1& 1 \\ 0& 1 \\ 0& 1 \end{pmatrix}	 .
\]
Then $\surjrep V \simeq \id$, given by the subspace $K e_1$ at each vertex: this is an epimorphic subrepresentation, and $V$ has no subrepresentations of dimension $(2, 1)$ or $(2, 2)$, so this subrepresentation must be maximal epimorphic.  Since each map is already injective, $\injrep V = V$, and so $\rkf V \simeq \id$.

The subrepresentation $\surjrep V$ is unique of dimension $(1, 1)$, but is not a direct summand: both uniqueness and the fact that $\surjrep$ has no complementary subrepresentation follow from the linear independence of the second columns of the above matrices.  So $\rkf V$ is not isomorphic to a direct summand of $V$.
In fact, there are no direct summands of any other dimension, and $V$ is actually indecomposable.
\end{example}

\section{Tensor Product and the Global Tensor Functor}
Since the tensor product of two surjective maps is surjective, and likewise for injective maps, the universal properties of $\surjrep$ and $\injrep$ induce natural transformations $\theta$ and $\zeta$ giving us commutative diagrams of functors
\begin{equation*}
\vcenter{\xymatrix{
&	V \otimes W \\
{\theta} : \surjrep V \otimes \surjrep W \ar@{^{`}->}[r] \ar@{^{`}->}[ur] & {\surjrep} (V \otimes W) \ar@{^{`}->}[u]} }
\qquad \qquad
\vcenter{\xymatrix{
V \otimes W \ar@{->>}[d] \ar@{->>}[dr]\\
{\zeta}: \injrep (V \otimes W) \ar@{->>}[r] & {\injrep} V \otimes \injrep W } } .
\end{equation*}
These natural transformations satisfy $D \circ \zeta = \theta \circ D$, and give $\surjrep$ and $\injrep$ the structure of weak tensor functors.
To show that $\surjrep$ is symmetric, that is, the structure map $\theta$ commutes with the commutativity constraint $\psi$ for $\surjrep$, one can verify that the following diagram (of natural transformations) is commutative:
\[
\vcenter{\xymatrix{
{\surjrep} V \otimes {\surjrep} W \ar@{^{`}->}[dd]^{\theta} \ar@{^{`}->}[dr] \ar[rrr]^{\psi} & & & {\surjrep} W \otimes {\surjrep} V \ar@{^{`}->}[dd]^{\theta} \ar@{^{`}->}[dl] \\
 & V \otimes W \ar[r]^{\psi} & W \otimes V \\
{\surjrep} (V \otimes W) \ar@{^{`}->}[ur] \ar[rrr]^{\surjrep(\psi)} & & & {\surjrep} (W \otimes V) \ar@{^{`}->}[lu]
}} .
\]
The left and right triangles commute because $\surjrep V \otimes \surjrep W$ is a subspace of $\surjrep (V \otimes W)$, by the universal property, so $\theta$ commutes with the monomorphisms of these functors to the identity functor on $\repq$.  The lower trapezoid commutes because $\surjrep$ is a subfunctor of $id_{\repq}$, and the upper trapezoid commutes because of the same statement, along with the fact that bifunctoriality of $\otimes$ forces $id_{\repq} \otimes id_{\repq} =  id_{\repq}$.  Checking the other conditions for $\surjrep$ and $\injrep$ to be weak tensor functors is similar.
However, the next example shows that neither $\surjrep$ nor $\injrep$ is a tensor functor.

\begin{example}
The maps $\theta$ and $\zeta$ are not in general isomorphisms.  For example, take $Q = \krontwomaps{\bullet}{\bullet}{}{}$ with representations 
\[
V = \krontwomaps{K}{K}{1}{0} \qquad \qquad W = \krontwomaps{K}{K}{0}{1} .
\]
Then $\surjrep V = \surjrep W = 0$, but $\surjrep (V \otimes W) = \krontwo{K}{0}$, so $\theta$ is not an isomorphism (dualize to get an analogous example for $\injrep$).
\end{example}

For linear maps $A$ and $B$, rank is multiplicative in the sense that $\rank (A \otimes B) = \rank A \cdot \rank B$.  Although we have just seen that neither $\surjrep$ nor $\injrep$ commutes with tensor product, the global tensor functor $\rkf$ does.

\begin{theorem}\label{rankmultthm}
There is a natural isomorphism of bifunctors
\begin{equation}\label{rankmulteq}
\rkf_Q V \otimes \rkf_Q W \cong \rkf_Q (V \otimes W)
\end{equation}
giving $\rkf_Q$ the structure of a tensor functor.
\end{theorem}

Before proving the theorem, we need to establish a technical lemma.  As usual, since $Q$ is fixed we omit this subscript.

\begin{lemma}
Consider the natural transformation of bifunctors defined by the composition
\begin{equation*}
\sigma: \surjrep( V \otimes W) \subseteq V \otimes W \onto V \otimes \injrep W
\end{equation*}
where the second map is $id_V \otimes q_W$, writing $q_W : W \to \injrep W$ for the canonical quotient.  Then 
\[
\im \sigma \subseteq \surjrep V \otimes \injrep W .
\]
\end{lemma}

\begin{proof}
We check this as maps of vector spaces at each vertex $z$.
Using the natural isomorphism of vector spaces
\[
\Hom_K (\surjrep(V \otimes W)_z,\, V_z \otimes_K \injrep (W)_z ) \cong \Hom_K (\surjrep(V \otimes W)_z \otimes_K \injrep (W)_z^*,\, V_z)
\]
we can identify $\sigma_z$ with the map 
\[
\pi_z : \surjrep(V \otimes W)_z \otimes \injrep (W)_z^* \to V_z \qquad \left(\sum v_i \otimes w_i \right) \otimes f \mapsto \sum f (w_i ) v_i
\]
which we want to show takes image in $\surjrep (V)_z$.

We claim that the subspace $M := \oplus_z \im \pi_z$ is an epimorphic subrepresentation of $V$.
To see this, let $a \in Q_1$ and set $x = ta,\ y = ha$.  Given 
\[
\sum f (w_i ) v_i = \pi_x \left[ \left( \sum v_i \otimes w_i \right) \otimes f \right] \in \im \pi_x
\]
where $\left(\sum v_i \otimes w_i \right) \in \surjrep(V \otimes W)_x$ and $f \in \injrep (W)_x^*$, we want to show that $V_{a}$ maps this element into $\im \pi_y$.  Now from (\ref{surjinjdualeq}), we have $\injrep (W)_x^* = \surjrep (DW )_x$, so there exists $g \in \surjrep( DW )_y$ such that $f = W_a^* g$.  Then 
\[
\begin{split}
V_a \left( \sum f (w_i ) v_i \right)&  = \sum f (w_i ) V_a (v_i) = \sum g\left( W_a (w_i) \right) V_a (v_i) \\
& = \pi_y \left[\left(\sum V_a (v_i) \otimes W_a (w_i) \right) \otimes g \right] \in \im \pi_y
\end{split}
\]
showing that $M$ is a subrepresentation of $V$.  

To see that this subrepresentation is epimorphic, a similar argument works.  Given
$$
\pi_y \left[\left(\sum s_j \otimes t_j \right) \otimes g \right] \in \im \pi_y
$$
where $\left(\sum s_j \otimes t_j \right) \in \surjrep(V \otimes W)_y$ and $g \in \injrep (W)_y^*$, there exists $\sum v_i \otimes w_i  \in \surjrep (V \otimes W)_x$ such that 
\[
(V_a \otimes W_a) \left(\sum v_i \otimes w_i \right) = \sum  V_a( v_i ) \otimes W_a (w_i) = \sum s_j \otimes t_j .
\]
Thus we have
\[
\begin{split}
V_a \left( \pi_x \left[\left(\sum  v_i \otimes w_i \right) \otimes W_a^* g \right] \right) & = V_a \left( \sum g\left(W_a (w_i)\right)\, v_i \right) = \sum g\left(W_a (w_i)\right) V_a (v_i ) \\
& = \pi_y \left[ \left(\sum V_a (v_i) \otimes W_a (w_i) \right) \otimes g \right] = \pi_y \left[ \left(\sum s_j \otimes t_j \right) \otimes g \right] .
\end{split}
\]
So we see that $V_a |_M $ is surjective for each arrow $a$, hence $M \subseteq \surjrep V$ by the universal property of $\surjrep$.
\end{proof}

\begin{proof}[Proof of Theorem \ref{rankmultthm}]
The lemma establishes that $\im \sigma \subseteq \surjrep V \otimes \injrep W$, giving the dashed arrow in the diagram
\begin{equation*}
\xymatrix{
{\surjrep} (V \otimes W ) \ar@{^{`}->}[d] \ar@{-->}[r]^{\sigma} & {\surjrep} V \otimes \injrep W \ar@{^{`}->}[d]  \ar@{->>}[r] & {\rkf V} \otimes \injrep W \ar@{^{`}->}[d] \\
V \otimes W \ar@{->>}[r] & V \otimes \injrep W \ar@{->>}[r] & {\injrep} V \otimes \injrep W } .
\end{equation*}
Here, every arrow represents a canonical natural transformation of bifunctors, but we will simply say ``map'' throughout the proof to avoid this cumbersome phrase.
Thus the map
\begin{equation*}
\surjrep (V \otimes W) \to \injrep V \otimes \injrep W
\end{equation*}
factors through $\rkf V \otimes \injrep W$.  But by applying the same reasoning it must also factor through ${\injrep} V \otimes \rkf W$, hence through the intersection (as subfunctors of $ \injrep V \otimes \injrep W$)
\begin{equation*}
({\injrep} V \otimes \rkf W) \cap ({\rkf} V \otimes \injrep W )= \rkf V \otimes \rkf W .
\end{equation*}

So we have a natural map
\begin{equation*}
\alpha : \surjrep (V \otimes W) \onto \rkf V \otimes \rkf W
\end{equation*}
which is surjective because the subrepresentation $\surjrep V \otimes \surjrep W$ already surjects onto the right hand side, by definition of the global tensor functor.  Applying the same argument with $Q^{op}$, then dualizing, we get a map
\begin{equation*}
\beta : \rkf V \otimes \rkf W \into \injrep (V \otimes W) .
\end{equation*}

We summarize this with the commutative diagram
\begin{equation*}
\xymatrix{
{\surjrep} V \otimes \surjrep W \ar@{->>} [r] \ar@{^{`}->} [d]_{\theta} & {\rkf} V \otimes \rkf W \ar@{^{`}->}[dr]^{\beta} \ar@{^{`}->} [r] & {\injrep} V \otimes \injrep W \\
{\surjrep} (V \otimes W) \ar@{->>}[ur]^{\alpha} \ar@{->>}[r] & {\rkf} (V \otimes W)  \ar@{^{`}->}[r] & {\injrep} (V \otimes W) \ar@{->>}[u]_{\zeta} }
\end{equation*}
which shows that the natural map
\begin{equation*}
{\surjrep} (V \otimes W) \into V \otimes W \onto {\injrep} (V \otimes W)
\end{equation*}
factors through $ {\rkf} V \otimes \rkf W$, and the universal property of the image of a map gives an isomorphism
\[
\rkf V \otimes \rkf W \cong \rkf (V \otimes W) .
\]
We already know that $\rkf (\id) \cong \id$, and to show that $\rkf$ satisfies the other conditions to be a tensor functor is straightforward.
\end{proof}

In particular, we have finally shown that in fact $\rankfunc_Q : \repq \to \rep(A_1)$ is a rank functor.

\section{Application of Rank Functions}\label{ranktensorsect}

For a given quiver $Q$, we can use the global rank functions of other quivers to construct rank functions on $Q$.  For example, if $P \subset Q$ is any connected subquiver, we can define a homomorphism from $R(Q) \to \Z$ by restricting $V$ to $P$ then applying the rank function of $P$.  This will be denoted by $r_P (V)$, with the restriction being understood.

More generally, for any map of directed graphs $\alpha : Q' \to Q$, we can push forward a rank function on $Q'$ to act on representations of $Q$.
In categorical language, such an $\alpha$ is nothing other than a functor $\alpha : \quivcat(Q') \to \quivcat(Q)$, and the pullback $\alpha^* V$ of a representation $V$ along $\alpha$ is just composition of functors
\[
\repfunc_{\alpha^* V} : \quivcat (Q') \xto{\alpha} \quivcat(Q) \xto{\repfunc_V} \vscat .
\]
In fact, $\alpha^* : \repq \to \rep(Q')$ is a tensor functor.  Since the composition of tensor functors is again a tensor functor, any rank functor $F$ on $Q'$ pushes forward to a rank functor on $Q$ by composition
\[
\alpha_* F := F \circ \alpha^* .
\]
In terms of representation rings, $\alpha$ induces a ring homomorphism $\alpha^* : R(Q) \to R(Q')$, so that when $f$ is a rank function on $\rep(Q')$, we get a pushforward rank function $(\alpha_* f) (V) := f(\alpha^* V)$ on $Q$.
The example of restriction above is just the case of the inclusion map $i: P \into Q$ of a subquiver, but this technique can be used to construct rank functions which don't come from any subquiver.  We illustrate with an example:

\begin{example}\label{coverqeg}
Consider the following two quivers, where the numbers at the vertices are just labels rather than dimension vectors:
\[
Q' = \QBmaps{1}{2}{3}{3}{4}{a}{b}{c}{c} \qquad \qquad Q = \QAmaps{1}{2}{3}{4}{a}{b}{c} .
\]
We have a map of directed graphs $\alpha : Q' \to Q$ by identifying vertices and arrows of the same label.  A representation
\[
V = \QAmaps{V_1}{V_2}{V_3}{V_4}{A}{B}{C} \in \repq
\]
pulls back to a representation 
\[
\alpha^* V = \QBmaps{V_1}{V_2}{V_3}{V_3}{V_4}{A}{B}{C}{C} \in \rep (Q')
\]
and it is apparent that $\alpha^*$ commutes with direct sum and tensor product.  Then we can compute the pushforward of the global rank function of $Q'$:
\[
\alpha_* r_{Q'} (V) = r_{Q'} (\alpha^* V) = \dim_K (\im CA \cap \im CB) .
\]
This function is distinct from the global rank function of $Q$, which can be computed from the definitions to be
\begin{equation*}
r_Q (V) = \rank \left(C |_{\im A \cap \im B} \right) .
\end{equation*}
Note that $\alpha_* r_{Q'}(V)= 0$ if $\supp V \neq Q_0$.

In Example \ref{tensororienteg2}, it was claimed that
we could find the decomposition of $W \otimes W$ in a different way than direct computation.  The multiplicativity of $r_Q$ and $\alpha_* r_{Q'}$ give us:
\begin{equation}\label{firsteq}
\alpha_* r_{Q'} \left( W \otimes W \right) =\alpha_* r_{Q'} \left( W \right)^2 = 1
\end{equation}
\begin{equation}\label{secondeq}
r_{Q} \left( W \otimes W \right) =r_{Q} \left( W \right)^2 = 0
\end{equation}
where the values of $\alpha_* r_{Q'} \left( W \right)$ and $r_Q \left( W \right)$ are computed using linear algebra with our description of $W$.
Then (\ref{firsteq}), along with additivity of $\alpha_* r_{Q'}$, implies that $W\otimes W$ has an indecomposable summand $Z$ such that $\supp Z = Q_0$.  But (\ref{secondeq}) implies that $\id_Q$ is not a direct summand of $W \otimes W$.  Since $\id_Q$ and $W$ are the only two indecomposable representations of $Q$ supported on all of $Q$, we get that $W$ must be a direct summand of $W \otimes W$, and the other indecomposable summands must be simple by dimension count.
\end{example}

We can think of this carefully chosen pushforward function as ``distinguishing'' $W$ from other indecomposable representations.  In the following examples, we show how to apply this idea to find structure in $R(Q)$.  These techniques can be used for many quivers (of finite, tame, or wild type) to study their representation rings.  We continue with examples of finite type in order to simplify the demonstrations.

\begin{example}\label{redrepringeg}
Continuing the previous example, there are 11 connected subquivers $P \subseteq Q$, each of which gives a rank function
\begin{equation*}
r_P : R(Q) \to \Z \qquad P \subseteq Q .
\end{equation*}
We also have $\alpha_* r_{Q'} : R(Q) \to \Z$.  These maps define a ring homomorphism to the product
\begin{equation*}
\Delta: R(Q) \to \Z^{12} .
\end{equation*}
This is actually an isomorphism, which is easy to check:  We already know that as an abelian group, $R(Q) \simeq \Z^{\oplus 12}$, with the indecomposable representations of $Q$ as a basis.  One can simply compute the value of $\Delta (V)$ for each indecomposable representation of $Q$, using explicit descriptions of the rank functors.  Then we verify that the image vectors form a $\Z$-module basis for $\Z^{12}$.

Thus the isomorphism class of a representation of $Q$ is completely determined by the values of these 12 explicitly given functions, and since this is a map of algebras, we can simplify the problem of computing tensor products of representations of $Q$ by just multiplying in $\Z^{12}$.
\end{example}

We cannot expect such simple structure for the representation rings of all quivers.  We saw with Examples \ref{tensororienteg1} and \ref{tensororienteg2} that the tensor product structure depends on orientation, and now we will see that the isomorphism class of $R(Q)$ depends on the orientation of $Q$.

\begin{example}
When $Q$ is the three subspace quiver of Example \ref{tensororienteg1}, we again have a surjective ring homomorphism
\begin{equation*}
\Gamma = \prod_{P \subseteq Q} r_P: R(Q) \to \Z^{11}
\end{equation*}
as in the previous example.  But one cannot find by inspection any other distinct rank function.  The analogous push-forward function for this orientation is actually equal to $r_Q$.  This might lead one to believe that $R(Q)$ has a non-reduced factor, which cannot be detected by homomorphisms to $\Z$.

Again, we can calculate the values $\Gamma (V)$ for each indecomposable $V$, and use linear algebra to see that $\Gamma$ is surjective with kernel generated by $E - F$, where
\[
E = \left( \threesubspaced{2}{1}{1}{1} \oplus \threesubspaced{1}{0}{0}{0} \right) \qquad F = \left(  \threesubspaced{1}{1}{0}{0} \oplus  \threesubspaced{1}{0}{1}{0} \oplus  \threesubspaced{1}{0}{0}{1} \right) .
\]
To compute powers of $E-F$, the only nontrivial multiplication we need is 
$\threesubspaced{2}{1}{1}{1} \otimes \threesubspaced{2}{1}{1}{1}$, which was done in Example \ref{tensororienteg1}.  Using this, it is easy to verify that $(E - F)^2 = 0$ in $R(Q)$, and
then that $\Gamma$ lifts to an isomorphism of $\Z$-algebras
\[
\tilde{\Gamma}: R(Q) \xto{\sim} \Z^{10} \times \frac{\Z[\varepsilon]}{\left(\varepsilon^2 \right)}
\]
where the map to the first factor is given by $\prod_{P \subsetneq Q} r_P$, and $\tilde{\Gamma} (E - F) = \varepsilon$.

\end{example}

So we see by example how rank functions can be used to give some basic structure of $R(Q)$, although there are not always enough tensor functors to completely determine $R(Q)$.

When $Q$ is not of finite representation type, one does not have the luxury of writing down a nice $\Z$-basis of $R(Q)$ on which we can easily compute the values of rank functions, so the above examples cannot be simply imitated.  However, one can use more complicated techniques with rank functions to study $R(Q)$.  For example, on certain quivers, one can combinatorially construct sets of ranks functions with additional structure, then show that this structure is reflected in $R(Q)$.

\bibliographystyle{plain}
\bibliography{ryanbiblio}

\def\cprime{$'$} \def\ocirc#1{\ifmmode\setbox0=\hbox{$#1$}\dimen0=\ht0
  \advance\dimen0 by1pt\rlap{\hbox to\wd0{\hss\raise\dimen0
  \hbox{\hskip.2em$\scriptscriptstyle\circ$}\hss}}#1\else {\accent"17 #1}\fi}
\begin{thebibliography}{10}

\bibitem{assemetal}
Ibrahim Assem, Daniel Simson, and Andrzej Skowro{\'n}ski.
\newblock {\em Elements of the representation theory of associative algebras.
  {V}ol. 1}, volume~65 of {\em London Mathematical Society Student Texts}.
\newblock Cambridge University Press, Cambridge, 2006.
\newblock Techniques of representation theory.

\bibitem{atiyahtall}
M.~F. Atiyah and D.~O. Tall.
\newblock Group representations, {$\lambda $}-rings and the
  {$J$}-homomorphism{$J$}-homomorphism.
\newblock {\em Topology}, 8:253--297, 1969.

\bibitem{BSriemannroch}
Armand Borel and Jean-Pierre Serre.
\newblock Le th\'eor\`eme de {R}iemann-{R}och.
\newblock {\em Bull. Soc. Math. France}, 86:97--136, 1958.

\bibitem{DMtannakian}
Pierre Deligne and James~S. Milne.
\newblock {\em Tannakian categories in Hodge cycles, motives, and {S}himura
  varieties}, volume 900 of {\em Lecture Notes in Mathematics}, pages 101--228.
\newblock Springer-Verlag, Berlin, 1982.
\newblock Available at \url{http://www.jmilne.org/math/articles/1982b.pdf}.

\bibitem{fultonreptheory}
William Fulton and Joe Harris.
\newblock {\em Representation theory}, volume 129 of {\em Graduate Texts in
  Mathematics}.
\newblock Springer-Verlag, New York, 1991.
\newblock A first course, Readings in Mathematics.

\bibitem{gabriel}
P.~Gabriel.
\newblock Unzerlegbare {D}arstellungen. {I}.
\newblock {\em Manuscripta Math.}, 6:71--103; correction, ibid. 6 (1972), 309,
  1972.

\bibitem{Gchernclasses}
Alexander Grothendieck.
\newblock La th\'eorie des classes de {C}hern.
\newblock {\em Bull. Soc. Math. France}, 86:137--154, 1958.

\bibitem{herschendaffinean}
Martin Herschend.
\newblock Solution to the {C}lebsch-{G}ordan problem for representations of
  quivers of type {$\tilde{\Bbb A}\sb n$}.
\newblock {\em J. Algebra Appl.}, 4(5):481--488, 2005.

\bibitem{herschtensorjpaa}
Martin Herschend.
\newblock Tensor products on quiver representations.
\newblock {\em J. Pure Appl. Algebra}, 212(2):452--469, 2008.

\bibitem{herschend07b}
Martin Herschend.
\newblock On the representation ring of a quiver.
\newblock {\em Algebr. Represent. Theory}, 2009.
\newblock to appear. Currently available online,
  \href{http://dx.doi.org/10.1007/s10468-008-9118-1}{\texttt{DOI:10.1007/s1046%
8-008-9118-1}}.

\bibitem{kacrootsystemsa}
V.~G. Kac.
\newblock Infinite root systems, representations of graphs and invariant
  theory.
\newblock {\em Invent. Math.}, 56(1):57--92, 1980.

\bibitem{knutsonlambda}
Donald Knutson.
\newblock {\em {$\lambda $}-rings and the representation theory of the
  symmetric group}.
\newblock Springer-Verlag, Berlin, 1973.
\newblock Lecture Notes in Mathematics, Vol. 308.

\bibitem{Mcategories}
Saunders Mac~Lane.
\newblock {\em Categories for the working mathematician}, volume~5 of {\em
  Graduate Texts in Mathematics}.
\newblock Springer-Verlag, New York, second edition, 1998.

\bibitem{popescuabeliancats}
N.~Popescu.
\newblock {\em Abelian categories with applications to rings and modules}.
\newblock Academic Press, London, 1973.
\newblock London Mathematical Society Monographs, No. 3.

\bibitem{segalequivktheory}
Graeme Segal.
\newblock Equivariant {$K$}-theory.
\newblock {\em Inst. Hautes \'Etudes Sci. Publ. Math.}, 34:129--151, 1968.

\bibitem{segallierepring}
Graeme Segal.
\newblock The representation ring of a compact {L}ie group.
\newblock {\em Inst. Hautes \'Etudes Sci. Publ. Math.}, 34:113--128, 1968.

\end{thebibliography}

\end{document}